\documentclass[12pt,a4paper,reqno]{amsart}
\usepackage{amssymb}
\usepackage{amscd}
\usepackage[pdftex,pdfpagelabels]{hyperref}
\usepackage{enumerate}
\usepackage{graphicx}
\usepackage{siunitx}
\usepackage{tikz-cd}
\usepackage{color}
\usepackage{MnSymbol}
\usepackage{url}
\hypersetup{colorlinks,linkcolor={black},citecolor={black},urlcolor={black}}  

\usetikzlibrary{arrows}
\numberwithin{equation}{section}

\usepackage{mathtools}
\usepackage[tableposition=top]{caption}
\usepackage{booktabs,dcolumn}

\DeclareFontFamily{OT1}{rsfs}{}
\DeclareFontShape{OT1}{rsfs}{n}{it}{<-> rsfs10}{}
\DeclareMathAlphabet{\mathscr}{OT1}{rsfs}{n}{it}

\theoremstyle{plain}

\newtheorem{theorem}{Theorem}[section]

\newtheorem{lemma}[theorem]{Lemma}
\newtheorem{corollary}[theorem]{Corollary}
\newtheorem{conjecture}[theorem]{Conjecture}

\theoremstyle{definition}

\newtheorem{definition}[theorem]{Definition}
\newtheorem{remark}[theorem]{Remark}

\newcommand\C{\mathbb{C}}
\newcommand\Q{\mathbb{Q}}

\begin{document}

\title[]{Effective upper bound of analytic torsion under Arakelov metric}

\author{Changwei Zhou}
\address{Binghamton University, Mathematics Department}
\email{czhou11@binghamton.edu}


\subjclass[2010]{14G40 , 58J52, 13D45}

\begin{abstract}  
Given a choice of metric on the Riemann surface, the regularized determinant of Laplacian (analytic torsion) is defined via the complex power of elliptic operators:
$$
\det(\Delta)=\exp(-\zeta'(0))
$$ 
In this paper we gave an asymptotic effective estimate of analytic torsion under Arakelov metric. In particular, after taking the logarithm it is asymptotically upper bounded by $g$ for $g>1$. 

The construction of a cohomology theory for arithmetic surfaces in Arakelov theory has long been an open problem. In particular, it is not known if $h^{1}(X,L)\ge 0$. We view this as an indirect piece of evidence that if such a cohomology theory exists, the $h^{1}$ term may be effectively estimated. 
\end{abstract}

\maketitle


\section{Introduction}

Let $(X_{\sigma},g)$ be a compact connected smooth Riemann surface without boundary (which we henceforth abbreviate as \emph{compact Riemann surface}). The metric Laplacian is defined to be
$$
\Delta_{g}(f)=\frac{1}{\sqrt{\det(g)}}\partial_{i}(\sqrt{\det(g)}g^{ij}\partial_{j}f)
$$
The regularized determinant of the metric Laplacian (which we henceforth abbreviate as \emph{analytic torsion}) is defined to be
$$
\det(\Delta_{g})=\exp(-\zeta'(0)), \zeta(s)=\frac{1}{\Gamma(s)}\int_{0}^{\infty}t^{s-1}\textrm{Tr}(e^{t\Delta}-P)dt
$$
where $P$ is the projection operator onto the kernel of the Laplacian.

The present paper presents the following theorem: 
\begin{theorem}
Let $\pi: X\rightarrow B=\textrm{Spec}(\mathcal{O}_{K})$ be an integral, flat, projective scheme of dimension $2$. Let $\sigma: K\rightarrow \C$ be a fixed archimedean place of $K$. Let $\Delta_{Ar}$ denotes the metric Laplacian associated to Arakelov metric on $X_{\sigma}$. 

We have the following effective estimate of the analytic torsion under Arakelov metric: For $g$ large enough ($g>10$, for example):
$$
-\infty<\log(\det(\Delta_{Ar}))<g
$$
In particular, for $g\ge 1$ the analytic torsion under Arakelov metric is always effectively upper bounded from above. The precise upper bound can be found in Appendix. 
\end{theorem}

Our proof of this result heavily used previous results by Jorgenson and Kramer, Wentworth and Wilms. The main idea behind the paper is that with Richard Wentworth and Jay Jorgenson's proof of the correct bosonization formula, we may write the scalar analytic torsion in terms of Faltings' delta function and area of the compact Riemann surface under Arakelov metric. With Robert Wilms's result, we may directly estimate Faltings' delta function. Thus to estimate the scalar analytic torsion, it is enough to bound the area of the surface. This was done essentially by Jorgenson and Kramer in their paper \cite{JK1}. This answers a question of Gillet and Soule in their paper \cite{Gillet-Soule2} in the setting of Arakelov metric. 

We note that it is crucial that we restrict ourselves to metric Laplacian under Arakelov metric. As we mentioned earlier, the paper \cite{Sarnak} showed such estimate is already very subtle even for hyperbolic metric. For example, the precise extreme metric for $g\ge 2$ is unknown. 

The organization of the paper is as follows: In section 2 we review the main background notation and previous work on this topic. In section 3-4 we present the proof and some ideas related to it, which is rather trivial. In Appendix section 5 we give an  that lists effective estimates for arithmetic surfaces of small genus except $g=1$. In Appendix section 6 we give a derivation for the case of elliptic curve. 

I thank A.Borisov, J. Bost, C. Soule, R. Wentworth, R. Wilms, J. Kramer, J. Jorgenson, P. Andreae, X.Wang, M. Kierlanczyk
 for generously sharing their ideas during the research process. In particular based Prof.Kramer's suggestion, I found the effective estimate. 

In May 30th, 2019 I received a preprint \cite{Jorgenson3} from Prof. Jorgenson. It is clear from the preprint that some of the results from this paper may be well known to Prof. Jorgenson and Prof. Wentworth as early as 1990. Further, his paper went much deeper than I did at here. 

Without the emotional support from friends like Jaiung Jun, Jingyu Zhao, Oleg Lazerev, Patrick Milano, Yiming Zhao and Zana Chan, the paper probably would never have been written down and I would be sweating on Leetcode instead. 

For obvious reasons, I am indebted to Prof. Vasiu. 

\section{Review of background material}
The main point for this section is to set notation and review work done by experts in the field. We first recall that Arakelov considered a metric on $X_{\sigma}$ whose area is 1:
\begin{definition}
Let $X_{\sigma}$ be a Riemann surface of genus $g\ge 1$. Let $\alpha_{j}$ be an orthonormal basis of holomorphic 1-forms on $X_{\sigma}$ with scalar product
$$
(\alpha,\beta)=\frac{i}{2}\int_{X}\alpha\wedge \overline{\beta}
$$
Then the volume form given by
$$
\mu_{\sigma}=\frac{i}{2g}\sum \alpha_{i}\wedge \overline{\alpha}_{j}
$$
is called the \it{canonical metric}.  
\end{definition}

We now proceed to define the Arakelov metric (\cite{Arakelov}). 
\begin{definition}\label{def}
Let 
$$
G_{\sigma}: X_{\sigma}\times X_{\sigma}\rightarrow \mathbb{R}
$$
be Green function if it satisfies the following properties:

\begin{itemize}
	\item 
	$G_{\sigma}$ is $C^{\infty}$ on $X_{\sigma}\times X_{\sigma}$. 
	\item 
	$G_{\sigma}(P,Q)$ has a first order zero on the diagonal, and it is non-zero if $P\not=Q$. 
	\item 
	If $P\not=Q$, then we require
	$$
	\partial_{P}\overline{\partial_{P}}\log G_{\sigma}(P,Q)=\pi i \mu_{\sigma}(P)
	$$
	\item 
	We require 
	$$
	\int_{X}\log G_{\sigma}(P,z)\mu_{\sigma}(z)=0
	$$
\end{itemize}
\end{definition}

Let $\sigma$ be an archimdean place of $K$. We fix the embedding $\sigma$ and let $M=X_{\sigma}$ be a fixed Riemann surface. Then we have so called Poincare residue isomorphism:
$$
\mathcal{O}_{M\times M}(-\Delta)\cong K_{M}
$$
where $\Delta$ is the diagonal. Using Green functions,  we may define a metric on the dual space $\mathcal{O}_{M\times M}(\Delta)$. Let $1_{\Delta}$ be its canonical section, we define the metric by
$$
(-\log |1_{\Delta}|)(P,Q)=\log G(P,Q)
$$
\begin{definition}\label{def}
The dual metric of above metric on $\mathcal{O}_{M\times M}(-\Delta)$ is the so called {Arakelov metric}. To be precise, we consider the diagonal morphism
$$
\mu: M\rightarrow M\times M
$$
and let $\mu^{*}$ be the pull-back metric on $\Omega^{1}(M)$. We want to mention the motivation of this construction (see \cite{Arakelov}): Classical adjunction formula suggests that 
$$
K_{D}=(K_{X}+D)|_{D}
$$
This can be seen from the conormal exact sequence induced from the inclusion map $i:D\rightarrow X$:
$$
0\rightarrow i^{*}\mathcal{O}(-D)\rightarrow i^{*}K_{X}\rightarrow K_{D}\rightarrow 0
$$
Arakelov wants the adjunction formula isomorphism to be an isometry for the case $D=P$:
$$
\Omega_{X}\otimes \mathcal{O}_{P}\cong \Omega_{P}
,
\eta\otimes \frac{s}{f}\rightarrow s\frac{\partial \eta}{\partial f}|_{f=0}
$$

Since the metric at $\mathcal{O}_{P}$ is given by $G(P,Q)$, we have the Arakelov metric on $\Omega_{X}$ to be
$$
|dz(P)|_{Ar}=\lim_{Q\rightarrow P}\frac{|z(Q)-z(P)|}{G(P,Q)}
$$
where $\eta=dz, s=1, f=z$. The Arakelov metric on $\mathcal{O}_{X}$ is defined to be its dual metric:
$$
\mu_{Ar}=\frac{|dz|}{|dz|_{Ar}}\frac{i}{2}dz\wedge \overline{dz}
$$
\end{definition}

In 1984, while working on an arithmetic analog of Riemann-Roch theorem for arithmetic surfaces, Faltings proved the arithmetic Noether formula \cite{Faltings}: As before we begin with a curve $C$ over a number field $K$. After a finite field extension, we may assume $C$ has semi-stable reduction over $K$. Further we denote by $\mathcal{X}$ the minimal regular model of $C$ over $B=Spec(O_{K})$. We have:
$$
12\widehat{\deg}\det \pi_{*}\omega_{\mathcal{X}/B}=\langle \omega_{X/B},\omega_{\mathcal{X}/B}\rangle+\sum_{v\in |B|}\delta_{v}\log (N(v))+\sum_{\sigma: K\rightarrow \mathbb{C}}\delta'(\mathcal{X}_{\sigma})
$$
By the geometric analogy, the interpretation of the formula is straightforward: The extra terms $\delta_{i}$ on the right hand side measures the degree of singularity of $\mathcal{X}$. Thus it should `blow up' on the boundary of the moduli space $M_{g}$. We note there is a discrepancy in normalization: The $\delta$-function we discuss at here, $\delta_{Fal}$ equals $\delta'+4g\log(2\pi)$. 

In 1987, Deligne proved a version of Riemann-Roch theorem using analytic torsion. In 1988 Bost, Gillet and Soule refined his result. They found that
$$
\delta_{Fal}(X)=-6D_{Ar}(X)+a(g_{X}), 
$$
where 
$$
D_{Ar}(X)=\log(\frac{\det(\Delta)_{Ar}}{\textrm{Area}_{Ar}(X)}),
$$
and
$$
a(g)=-8g\log(2\pi)+(1-g)(-24\zeta'_{\mathbb{Q}}(-1)+1-6\log(2\pi)-2\log(2))
$$
It is thus clear that $\delta_{Fal}(X)$'s growth only depend on the analytic torsion and the area of $X$ under Arakelov metric. The original computation of $a(g)$ in Soule's paper \cite{Soule2} turned out to be incorrect, most notably an sign error $(g-1)$ instead of $1-g$. Later in 2008, Richard Wentworth re-derived the correct formula rigorously in his paper \cite{Wentworth}. We wish to mention there is an independent proof by Jay Jorgenson \cite{Jorgenson3}.

In \cite{JK2001} and \cite{JK1}, Jorgenson and Kramer systematically investigated the area of a compact Riemann surface under Arakelov metric. We note that everything below is only valid for $g>1$. Their result, presented in compact form (page 7, Corollary 3.3), is that
$$
\log(\frac{\mu_{Ar}(z)}{\mu_{hyp}(z)})\le 4\pi(1-\frac{1}{g})\int^{\infty}_{1}K_{H}(t,0)dt-\frac{c_{sel}}{g(g-1)}+\frac{1}{g-1}-\log[4]
$$
where $\mu_{hyp}(z)$ denotes the hyperbolic metric and $\mu_{Ar}(z)$ denotes the Arakelov metric. Here the constant $K_{H}(t,0)$ is the hyperbolic heat kernel.

In 2016, Robert Wilms proved in his paper (\cite{Wilms}) that 
$$
\delta_{Fal}(X)>-2g\log(2\pi^4)
$$
for all $g$. We refer the reader to his thesis for a detailed exposition of his paper. 

\section{A heuristic interpretation of the result}
Before we elaborate the proof, we want to offer a few comments why the theorem 'should be expected' other than that it fits into arithmetic intersection theory framework established by Faltings and Gillet-Soule. In fact, Jorgenson and Kramer's earlier result in \cite{JK2} indicates that the area of the surfaces shrinks to $0$ near the boundary of $M_{g}$. Thus an easy argument using the compactness of $M_{g}$ showed the upper bound exists and only dependent on $g$ for $g>1$. The real surprise is that the upper bound is only linear. This feels rather mysterious at first sight. 

Here we note that using gluing formula of determinant of elliptic operators we have
$$
\frac{\det_{Ar}(X)}{Area_{Ar}(X)}=\det(X_{+}, Ar)\det(X_{-}, Ar)\frac{\det(N_{l,\gamma, h})}{\textrm{length}_{Ar}(\gamma)}
$$
The notation is borrowed from Wentworth's paper: $\gamma$ is a simple separating curve on $X$ that separates $X$ into two parts, $X_{+}$ and $X_{-}$. The determinant on the right hand side is evaluated in terms of Laplacian with Dirichlet boundary condition on $\gamma$. The symbol $N_{l,\gamma}$ denotes the Newmann jumping operator. 

An `obvious idea' now is to take logarithm, and use induction on $g$ to control the growth of each term. However, the Arakelov metric inherited this way for the surfaces with boundary of genus $g-1$ is not the same Arakelov metric for closed surface $g$. So far more involved analysis similar to Wentworth's work is necessary. We conjecture that the term involving the Newmann jumping operator is bounded from above. Heuristically, the determinant `factorizes' and as $g$ grows the extra contribution should be at most linear. 

\section{Proof of main theorem}
The associated integral may be explicitly bounded by
$$
\int^{\infty}_{1}K_{H}(t,0)dt\le \int^{\infty}_{1}\frac{e^{-t/4}}{4\pi t}dt\le 0.0832<0.1
$$
Thus for $g>1$, the constant term in the above sum may be evaluated to be bounded above by $\frac{1}{g-1}-0.33\le 0.67\le 1$. The only term left to estimate is the $c_\text{sel}$ for hyperbolic surfaces, which Jorgenson and Kramer estimated in \cite{JK2001} to be (page 13, remark 3.5):
$$
c_\text{sel}\ge -4\log(1366(g-1))
$$
As a result we have
$$
\log(\frac{\mu_{Ar}(z)}{\mu_{hyp}(z)})\le 1+\frac{4\log(1366(g-1))}{g(g-1)}
$$
Exponentiate both sides and integrate, we have
$$
\textrm{Area}(X_{\sigma})< e*4\pi(g-1)* (1366(g-1))^{\frac{4}{g(g-1)}}< 36 (g-1)(1366(g-1))^{\frac{4}{g(g-1)}}
$$
where we used the fact for hyperbolic surfaces the area is given by $4\pi(g-1)$ (Gauss-Bonnet). 

By combing Wilms' result and Wentworth's result mentioned earlier, we have
$$
-6D_{Ar}(X)+a(g_{X})>-2g\log(2\pi^4)
$$
where 
$$
D_{Ar}(X)=\log(\frac{\det(\Delta)_{Ar}}{\textrm{Area}_{Ar}(X)}),
$$
and
$$
a(g)=-8g\log(2\pi)+(1-g)(-24\zeta'_{\mathbb{Q}}(-1)+1-6\log(2\pi)-2\log(2))
$$
To ease notation we ignore the $Ar$ in the subscript. Rearranging and clearing off the terms, we have
$$
\log(\det(\Delta))<\frac{\log(2\pi^4)}{3}g+\frac{1}{6}a(g)+\log(Area(X))
$$
By previous computation we have
$$
\log(Area(X))<\log[36]+\log[g-1]+\frac{4}{g(g-1)}\log(1366(g-1))\sim O(\log(g))
$$
We now focus on the leading term of size $O(g)$. Combining with the first term we have an effective asymptotic upper bound to be
$$
(\frac{\log(2\pi^4)}{3}-8/6*\log[2\pi]-\frac{1}{6}(-24\zeta'_{\mathbb{Q}}(-1)+1-6\log(2\pi)-2\log(2)))g
$$
The constants can be explicitly evaluated to be
$$
\approx 1.7573-2.4505+2.07-\frac{1}{6}+4\zeta_{\mathbb{Q}}'(-1)< 1.21-0.67=0.56<1
$$
As a result we have the final effective estimate to be (valid for all $g>1$:)
$$
-\infty < \log(\det(\Delta))< 0.56g+E(g)
$$
where $E(g)\le 2\log(g)<0.44g$ for $g$ large enough. As a result the whole term is asymptotically bounded by $g$. We have the following effective formula for $E(g)$: $$E(g)=\log[36]+\log[g-1]+\frac{4}{g(g-1)}\log(1366(g-1))+\frac{1}{6}(-24\zeta'_{\mathbb{Q}}(-1)+1-6\log(2\pi)-2\log(2))$$
and a more refined estimate can be given by
$$
E(g)= \frac{1}{g-1}+\log(g-1)+\frac{4}{g(g-1)}\log(1366(g-1))+\frac{1}{6}(-24\zeta'_{\mathbb{Q}}(-1)+1-6\log(2\pi)-2\log(2))+2.1890125
$$
In particular arithmetic computation showed that for $g>10$, we have $E(g)<0.44g$. Thus for the above estimate, it suffice to let $g>10$.

\begin{corollary}
For a compact Riemann surface of genus $g\ge 1$, the difference of Faltings metric and Quillen metric's logarithm has an asymptotic lower bound by a constant. The constant only depends on $g$. 
\end{corollary}
\proof
By Jorgenson's result in his preprint (page 40) we have
$$
h_{F}(L)-h_{Q}(L)=\log(Area)-\log(\det(\Delta))-2\log(2\pi)
$$
where $h_{F}(L), h_{Q}(L)$ denotes the logarithm of Faltings and Quillen metric respectively. Re-write it we have
\begin{align*}
h_{F}(L)-h_{Q}(L)&\ge -Dr_{Ar}(X)-2\log(2\pi)\\
&= -2 \log(2\pi)-\frac{\alpha(g)}{6}-\frac{g}{3}\log(2\pi^4)\\
&=-2\log(2\pi)-\frac{-8g\log(2\pi)+(1-g)(-24\zeta'_{\mathbb{Q}}(-1)+1-6\log(2\pi)-2\log(2))}{6}-\frac{g}{3}\log(2\pi^4)\\
&=(\frac{4}{3}\log(2\pi)-\log(2\pi^4)/3+4\zeta'_{\Q}(-1)-\frac{1}{6}+\log(2\pi)+\frac{1}{3}\log(2))g+C
\end{align*}
where $C$ is the effective constant
$$
C=-\log(2\pi)+4\zeta'_{\mathbb{Q}}(-1)-\frac{1}{6}+\frac{1}{3}\log(2)\approx -3.6113717392987086-0.661685\ge -4.273
$$
And the coefficient for $g$ is bounded below by
$$
(\frac{4}{3}\log(2\pi)-\log(2\pi^4)/3+4\zeta'_{\Q}(-1)-\frac{1}{6}+\log(2\pi)+\frac{1}{3}\log(2))\approx 1.933721640489272\ge 1.934
$$
As a result the difference can be estimated to be
$$
h_{F}(L)-h_{Q}(L)\ge 1.934g-4.273> -2.334
$$
effective for all $g\ge 1$. 
\qed 

\remark 
We note that obtaining an effective upper bound on $M_{g}$ would be much more difficult as we would have to study the degeneration of $\delta_{X}$ in $M_{g}$. For recent work on this we cite \cite{de Jong2}. 

\remark 
By Elkies' result, the Faltings metric is 'much larger' than the $L^{2}$ metric when $\deg(L)$ is large enough. It would be interesting to know if semi-positive condition of the curvature under Arakelov metric induces an effective estimate of analytic torsion.

\section{Appendix: Effective bounds for small genus}
For the convenience of the reader we list some of the upper bounds we obtained for small $g$. We note that the case for $g=0$ was slightly misleading as Arakelov metric in this case is not exactly the round metric of the unit sphere (it is a multiple of Fubini-Study metric), which caused the error in Soule's Bourbaki paper. We did the numerical computation using Wolfram Alpha:

\begin{itemize}
\item $g=0$: Computational result by \cite{Jorgenson3} is
$$
\det(\Delta_{g})=\exp(-4\zeta'(-1)+7/6-4/3\log[2])\approx 2.46984
$$
\item 
$g=1$: After identifying $X_{\sigma}=\langle 1, \tau\rangle $ with $\tau=\langle x, iy\rangle$ we have: $$\det(\Delta_{g})\le 2\pi y^2 e^{-\frac{\pi}{2} y+\frac{3}{\pi y}}$$
(For an elementary derivation, see section 6)
\item $g=2$: 18.01100181
\item $g=3$: 9.76363
\item $g=4$: 8.13548
\item $g=5$: 7.88854
\item $g=6$: 8.10036
\item $g=7$: 8.50296
\item $g=8$: 8.99605
\item $g=9$: 9.53577
\item $g=10$: 10.1007
\item $10\le g\le 3580$: Bounded above by $g$. 
\item $g\ge 3580$: Bounded above by $0.5474277074g+1$. 
\end{itemize}

\section{Appendix: The case for elliptic curve}
We identify the elliptic curve to be generated by the lattice $\langle 1, \tau\rangle$:

\lemma
Let $\tau=x+iy$. Then we have the Arakelov metric associated to the elliptic curve to be
$$
|dz|_{Ar}=\frac{\sqrt{y}}{2\pi}||\eta(\tau)||^{-2}=\frac{\sqrt{y}}{2\pi}*\frac{1}{\sqrt{y}}\frac{1}{\eta(\tau)^2}=\frac{1}{2\pi \eta(\tau)^2}
$$
\proof 
This was proved by Faltings in \cite{Faltings}.  
\qed

\lemma 
The Arakelov metric over functions is given by
$$
\mu_{Ar}=\pi|\eta(\tau)|^2 dz\wedge \overline{dz}
$$
\proof
Since the metric on $\mathcal{O}_{X}$ is the dual metric of the metric on $K_{X}$, we have
$$
\mu_{Ar}=\frac{|dz|}{|dz|_{Ar}}\frac{i}{2}dz\wedge \overline{dz}=\pi i|\eta(\tau)|^2dz\wedge \overline{dz}
$$

\lemma
The Arakelov area associated to the Arakelov metric is
$$
\textrm{Area}_{Ar}=2\pi y |\eta(\tau)|^2
$$
\proof
We know that the area associated to the fundamental domain is
$$
\frac{i}{2}\int_{D}dz\wedge \overline{dz}
$$
As a result, we have
$$
\textrm{Area}_{D}=\int^{x+iy}_{0}\int^{1}_{0}|\mu|_{Ar}=2\pi y |\eta(\tau)|^2
$$
\qed 

\lemma
If we normalize the unit area to be $1$, such that $E\cong \mathbb{R}^{2}/L$, where $L=\mu\langle 1, z\rangle$, then the regularized determinant associated to the flat metric for a lattice of unit area is
$$
\det(\Delta)=y|\eta(z)|^{4}
$$
\proof
The computation here was done by Osgood, Philips and Sarnark. They pointed out that for flat metric with area $1$ on the complex plane, we have (\cite{Sarnak}, page 204):
$$
\log(\det(\Delta))=\log(y|\eta(z)|^4)
$$
\qed 

\lemma 
The regularized determinant associated to the Arakelov metric is
$$
\log (\det(\Delta_{Ar}))=\log(2\pi)+2\log(y)+6\log(|\eta(\tau)|)
$$
\proof
By Lemma 6.2, the Arakelov metric is a scaled version of the flat metric on $\langle 1, \tau\rangle$. By Lemma 6.3, the scaling factor versus unit area metric is 
$$
\gamma=\sqrt{2\pi y |\eta(\tau)|^2}
$$
By Lemma 6.4, the regularized determinant associated to the lattice with unit area flat metric is
$$
\log(\det({\textrm{Unit}}))=\log(y|\eta(\tau)^4)|)
$$
We now cite a result in (\cite{Sarnak}, page 156): For elliptic curve, under a scaling factor of $\gamma^2: g(x,x)\rightarrow \gamma^2g(x,x)$ of the metric $g$, we have
$$
\log(\det(\gamma^2 g))=2\log(\gamma)+\log(\det(g))
$$
Indeed this can be proved directly: We have (note $Z(0)=-1$):
$$
Z_{L}(s)=\sum_{l\in L^{*}}\frac{1}{(2\pi |l|)^{2s}}\rightarrow Z_{\gamma L}=\gamma^{2s}\sum_{l\in L^{*}}\frac{1}{(2\pi |l|)^{2s}}
$$
Combining everything together we have
$$
\log (\det(\Delta_{Ar}))=\log(2\pi y |\eta(\tau)|^2)+\log(y|\eta(\tau)^4|)=\log(2\pi)+2\log(y)+6\log(|\eta(\tau)|)
$$
\qed 

\corollary
We have
$$
\log(y^2|\eta(\tau)|^6)=\log (y^2 e^{-\frac{\pi}{2} y}\prod^{\infty}_{n=1}|1-e^{2\pi i n\tau}|^6)
$$
to be bounded above by 
$$
2\log(y)-\frac{\pi y}{2}+\frac{6}{\pi y}
$$
\proof
We have the following elementary inequality (I learned this by reading Rafael von Kanel):
$$
\log|(\prod^{\infty}_{n=1}(1-q^{n}))|\le \frac{|q|}{(1-|q|)}, q=e^{2\pi i \tau} 
$$
As the result the above term is bounded above by
$$
2\log(y)-\frac{\pi y}{2}+6\frac{e^{-2\pi y}}{1-e^{-2\pi y}}= 2\log(y)-\frac{\pi y}{2}+\frac{6}{e^{2\pi y}-1}\le 2\log(y)-\frac{\pi y}{2}+\frac{3}{\pi y}
$$
where we recall $q=e^{2\pi i \tau}=e^{2\pi i x}*e^{-2\pi  y}$, so we have $|q|=e^{-2\pi y}$.  
\qed

\end{document}